# Disturbance decoupled functional observers for fault estimation in nonlinear systems*

Sunjeev Venkateswaran, Costas Kravaris*

*Abstract*— This work deals with the problem of designing disturbance decupled observers for the estimation of a function of the states in nonlinear systems. Necessary and sufficient conditions for the existence of lower order disturbance decoupled functional observers with linear dynamics and linear output map are derived. Based on this methodology, a fault-estimation scheme based on disturbance decoupled observers will be presented. Throughout the paper, the application of the results will be illustrated through a chemical reactor case study.

## I. INTRODUCTION

In control theory, a functional observer is an auxiliary system that is driven by the available measured outputs and mirrors the dynamics of a physical process, in order to estimate one or more functions of the system states [1, 2]. Besides being of theoretical importance, the use of functional observers arises in many applications. For example, functional estimates are useful in feedback control system design because the control signal is often a function of the states, and it is possible to utilize a functional observer to directly estimate the feedback control signal[1-3].

Over the past fifty years, considerable research effort has been carried out on estimating functions of the state vector for linear systems ever since Luenberger introduced the concept of functional observers in 1966[2] and proved that it is feasible to construct a functional observer with number of states equal to observability index minus one. Subsequent research has focused on lower order functional observers where necessary and sufficient conditions for their existence and stability have been derived[4-6], and parametric approaches to the design of lower order functional observers[7] and algorithms for solving the functional observer design conditions have also been developed[4, 5, 8] For nonlinear systems, functional observers for Lipschitz systems [9, 10] and a class of nonlinear systems that can be decomposed as sum of Lipschitz and non-Lipschitz parts [7] (with the non-Lipschitz part considered as an unknown input/disturbance) have been developed. More recently, the problem of designing functional observers for estimating a single nonlinear functional has been tackled for general nonlinear systems from the point of view of observer error linearization [3, 11] and the approach has been extended to a disturbance decoupled fault detection and isolation [12].

An important issue in the design of functional observers arises from the fact that accurate modeling of a real system is difficult and unknown disturbances/ uncertainties could hamper convergence of the estimate to its true value. This motivates imposing disturbance decoupling requirements in the design of a functional observer.

The present work extends the results presented in [3, 11] by considering the effect of unknown external disturbances on functional observer design in nonlinear systems. Necessary and sufficient conditions for disturbance decoupling are derived for functional observers designed from the point of view of observer error linearization. Following this, the method will be applied to fault estimation, where it will be assumed that the disturbances follow the dynamics of an exo-system. This is in the same vein as disturbance generators, in the context of the internal model principle, for regulation problems to achieve robust asymptotic regulation[13-15]. In our case, knowledge of how the disturbances is generated will enable robust tracking the fault signal.

Section 2 will review the observer error linearization problem studied in [3, 11]. Section 3 will study the effect of disturbance decoupling on functional observer design and necessary and sufficient conditions for disturbance decoupling will be derived. The disturbance decoupled functional observer design methodology will be illustrated via a non-isothermal Continuous Stirred Tank Reactor (CSTR) case study. Following this, the approach will extended to fault estimation and will be applied to the CSTR case study.

## II. FUNCTIONAL OBSERVERS WITH LINEAR ERROR DYNAMICS IN THE ABSENCE OF DISTURBANCES

In this section, a review of the design of functional observers from the point of view of observer error linearization will be given. Necessary and sufficient conditions for the existence of functional observers with linear error dynamics derived in [11] are presented, leading to simple formulas for observer design with eigenvalue assignment.

*Functional Observer Linearization Problem*
Given a system of the form

$$\frac{dx}{dt} = F(x) \quad (2.1)$$

$$y = H(x)$$
$$z = q(x)$$

where $F: \mathbb{R}^n \to \mathbb{R}^n, H: \mathbb{R}^n \to \mathbb{R}^p, q: \mathbb{R}^n \to \mathbb{R}$ are smooth nonlinear functions, y is the vector of measured outputs and z is the output to be estimated, find a functional observer of the form

*Research supported by National Science Foundation through the grant CBET-2133810. Both authors are with Texas A&M University, College Station, TX, 77843, USA (corresponding author e-mail: kravaris@tamu.edu).

$$\frac{d\hat{\xi}}{dt} = A\hat{\xi} + By \qquad (2.2)$$
$$\hat{z} = C\hat{\xi} + Dy$$

where A, B, C, D are $\nu \times \nu, \nu \times p, 1 \times \nu, 1 \times p$ matrices respectively, with A having stable eigenvalues. Equivalently, find a continuously differentiable mapping $T: \mathbb{R}^n \to \mathbb{R}^\nu$ such that

$$\frac{\partial T}{\partial x}(x)F(x) = AT(x) + BH(x) \qquad (2.3)$$

and

$$q(x) = CT(x) + DH(x) \qquad (2.4)$$

Assuming that the above problem can be solved, the resulting error dynamics will be linear:

$$\frac{d}{dt}\left(\hat{\xi} - T(x)\right) = A\left(\hat{\xi} - T(x)\right) \qquad (2.5)$$
$$\hat{z} - z = C(\hat{\xi} - T(x))$$

from which $\hat{z}(t) - z(t) = Ce^{At}(\hat{\xi}(0) - T(x(0)))$. With the matrix A having stable eigenvalues, the effect of the initialization error $\hat{\xi}(0) - T(x(0))$ will die out, and $\hat{z}(t)$ will approach $z(t)$ asymptotically.

The main result from previous paper [3, 11] is outlined in the following proposition:

Proposition 1: *For a nonlinear system of the form (2.1), there exists a functional observer of the form*

$$\frac{d\hat{\xi}}{dt} = A\hat{\xi} + By \qquad (2.2)$$
$$\hat{z} = C\hat{\xi} + Dy$$

*with the eigenvalues of A being the roots of a given polynomial* $\lambda^\nu + \alpha_1 \lambda^{\nu-1} + \cdots + \alpha_{\nu-1}\lambda + \alpha_\nu$, *if and only if*
$L_F^\nu q(x) + \alpha_1 L_F^{\nu-1} q(x) + \cdots + \alpha_{\nu-1} L_F q(x) + \alpha_\nu q(x)$
*is $\mathbb{R}$-linear combination of*
$H_j(x), L_F H_j(x), \ldots, L_F^\nu H_j(x), j = 1, \cdots, p$, *where*
$L_F = \sum_{k=1}^n F_k(x) \frac{\partial}{\partial x_k}$ *denotes the Lie derivative operator.*

In particular, if there exist constant row vectors $\beta_0, \beta_1, \cdots, \beta_\nu \in \mathbb{R}^p$ such that
$L_F^\nu q(x) + \alpha_1 L_F^{\nu-1} q(x) + \cdots + \alpha_{\nu-1} L_F q(x) + \alpha_\nu q(x)$
$= L_F^\nu(\beta_0 H(x)) + L_F^{\nu-1}(\beta_1 H(x)) + \cdots + L_F(\beta_{\nu-1} H(x))$
$+ \beta_\nu H(x)$
$= \sum_{j=1}^p \left(\beta_{0_j} L_F^\nu H_j(x) + \beta_{1_j} L_F^{\nu-1} H_j(x) + \cdots + \beta_{(\nu-1)_j} L_F H_j(x) + \beta_{\nu_j} H_j(x)\right)$ (2.6)

the mapping
$$T(x) = \begin{bmatrix} \left(\begin{array}{c} L_F^{\nu-1} q(x) + \alpha_1 L_F^{\nu-2} q(x) + \cdots + \alpha_{\nu-1} q(x) \ldots \\ -L_F^{\nu-1}(\beta_0 H(x)) - L_F^{\nu-2}(\beta_1 H(x)) - \cdots - \beta_{\nu-1} H(x) \end{array}\right) \\ \vdots \\ L_F q(x) + \alpha_1 q(x) - L_F(\beta_0 H(x)) - \beta_1 H(x) \\ q(x) - \beta_0 H(x) \end{bmatrix} \qquad (2.7)$$

satisfies conditions (2.3) and (2.4) with

$$A = \begin{bmatrix} 0 & 0 & \cdots & 0 & -\alpha_\nu \\ 1 & 0 & \cdots & 0 & -\alpha_{\nu-1} \\ \vdots & \ddots & \ddots & \vdots & \vdots \\ 0 & \cdots & 1 & 0 & -\alpha_2 \\ 0 & \cdots & 0 & 1 & -\alpha_1 \end{bmatrix}, \quad B = \begin{bmatrix} \beta_\nu - \alpha_\nu \beta_0 \\ \beta_{\nu-1} - \alpha_{\nu-1}\beta_0 \\ \beta_{\nu-2} - \alpha_{\nu-2}\beta_0 \\ \vdots \\ \beta_1 - \alpha_1 \beta_0 \end{bmatrix}$$
$$C = [0\ 0 \cdots 0\ 1], \qquad D = \beta_0 \qquad (2.8)$$

Hence, with the above A, B, C and D matrices, system (2.2) is a functional observer for (2.1).

### III. DISTURBANCE DECOUPLED FUNCTIONAL OBSERVERS WITH LINEAR ERROR DYNAMICS

Consider a nonlinear system that is subject to disturbance inputs $W \in \mathbb{R}^m$

$$\frac{dx}{dt} = F(x) + E(x)W \qquad (3.1)$$
$$y = H(x) + K(x)W$$
$$z = q(x)$$

where $F(x), H(x), q(x)$ $E(x)$ and $K(x)$ are smooth nonlinear functions. The objective is to build a functional observer whose output $\hat{z}$ is unaffected by the inputs W. In the context of the functional observer linearization problem, we seek for a functional observer of the form (2.2) with $T(x)$ solution of (2.3) and (2.4). Then, the resulting error dynamics will satisfy

$$\frac{d}{dt}\left(\hat{\xi} - T(x)\right) = A\left(\hat{\xi} - T(x)\right)$$
$$+ \left(BK(x) - \frac{\partial T}{\partial x}(x)E(x)\right)W \qquad (3.2)$$
$$\hat{z} - q(x) = C\left(\hat{\xi} - T(x)\right) + DK(x)W$$

The functional observer error will be completely unaffected by the disturbance inputs W if the coefficients of W in (3.2) vanish:

$$-\frac{\partial T}{\partial x}(x)E(x) + BK(x) = 0 \qquad (3.3)$$
$$DK(x) = 0 \qquad (3.4)$$

A functional observer of the form (2.2) that satisfies the above conditions will be called *disturbance decoupled* functional observer.

Proposition 2: *For a nonlinear system of the form (3.1), there exists a disturbance decoupled functional observer of the form (2.2) with the eigenvalues of A being the roots of a given polynomial* $\lambda^\nu + \alpha_1 \lambda^{\nu-1} + \cdots + \alpha_{\nu-1}\lambda + \alpha_\nu$, *if and only if there exist constant row vectors* $\beta_0, \beta_1, \cdots, \beta_\nu \in \mathbb{R}^p$, *that satisfy the following conditions:*

a) $L_F^\nu q(x) + \alpha_1 L_F^{\nu-1} q(x) + \cdots + \alpha_{\nu-1} L_F q(x) + \alpha_\nu q(x)$
   $= L_F^\nu(\beta_0 H(x)) + L_F^{\nu-1}(\beta_1 H(x)) + \cdots$
   $+ L_F(\beta_{\nu-1} H(x)) + \beta_\nu H(x)$

b) $\begin{cases} \sum_{\ell=0}^{\nu-\kappa} L_E L_F^{\nu-\kappa-\ell}(\beta_\ell H(x)) + \beta_{\nu-\kappa+1} K(x) = \\ L_E L_F^{\nu-\kappa} q(x) + \sum_{\ell=1}^{\nu-\kappa} \alpha_\ell L_E L_F^{\nu-\kappa-\ell} q(x), \kappa = 1, \cdots, \nu \\ \beta_0 K_i(x) = 0 \end{cases}$

Proof: Condition (a) is exactly the condition of Proposition 1 (equation (2.6)), which is necessary and sufficient for (2.2) to be a functional observer for (3.1) in the absence of disturbances (W=0). In order to be a disturbance decoupled functional observer, conditions (3.3) and (3.4) need to be applied. These can be written component-wise as follows:

$$\frac{\partial T_1(x)}{\partial x}E(x) - B_1 K(x) = 0$$
$$\frac{\partial T_2(x)}{\partial x}E(x) - B_2 K(x) = 0$$
$$\vdots$$
$$\frac{\partial T_{\nu-1}(x)}{\partial x}E(x) - B_{\nu-1} K(x) = 0$$
$$\frac{\partial T_\nu(x)}{\partial x}E(x) - B_\nu K(x) = 0$$
$$D K(x) = 0$$

Substituting the expressions for B, D and T(x) from (2.8) and (2.7) to the above equations lead to the following conditions:

$$L_E L_F^{\nu-1}(\beta_0 H(x)) + \cdots + L_E L_F(\beta_{\nu-2} H(x)) + L_E(\beta_{\nu-1} H(x)) + \alpha_\nu \beta_0 K(x) - \beta_\nu K(x)$$
$$= L_E L_F^{\nu-1} q(x) + \alpha_1 L_E L_F^{\nu-2} q(x) + \cdots + \alpha_{\nu-1} L_E q(x)$$
$$L_E L_F^{\nu-2}(\beta_0 H(x)) + \cdots + L_E L_F(\beta_{\nu-3} H(x)) + L_E(\beta_{\nu-2} H(x)) + \alpha_{\nu-1} \beta_0 K(x) - \beta_{\nu-1} K(x)$$
$$= L_E L_F^{\nu-2} q(x) + \alpha_1 L_E L_F^{\nu-3} q(x) + \cdots + \alpha_{\nu-2} L_E q(x)$$
$$\vdots$$
$$L_E L_F(\beta_0 H(x)) + L_E(\beta_1 H(x)) + \alpha_2 \beta_0 K(x) - \beta_2 K(x)$$
$$= L_E L_F q(x) + \alpha_1 L_E q(x)$$
$$L_E(\beta_0 H(x)) + \alpha_1 \beta_0 K(x) - \beta_1 K(x) = L_E q(x)$$
$$\beta_0 K(x) = 0$$

which can be written equivalently as

$$L_E L_F^{\nu-1}(\beta_0 H(x)) + \cdots + L_E L_F(\beta_{\nu-2} H(x)) + L_E(\beta_{\nu-1} H(x)) - \beta_\nu K(x)$$
$$= L_E L_F^{\nu-1} q(x) + \alpha_1 L_E L_F^{\nu-2} q(x) + \cdots + \alpha_{\nu-1} L_E q(x)$$
$$L_E L_F^{\nu-2}(\beta_0 H(x)) + \cdots + L_E L_F(\beta_{\nu-3} H(x)) + L_E(\beta_{\nu-2} H(x)) - \beta_{\nu-1} K(x)$$
$$= L_E L_F^{\nu-2} q(x) + \alpha_1 L_E L_F^{\nu-3} q(x) + \cdots + \alpha_{\nu-2} L_E q(x)$$
$$\vdots$$
$$L_E L_F(\beta_0 H(x)) + L_E(\beta_1 H(x)) - \beta_2 K(x)$$
$$= L_E L_F q(x) + \alpha_1 L_E q(x)$$
$$L_E(\beta_0 H(x)) - \beta_1 K(x) = L_E q(x)$$
$$\beta_0 K(x) = 0$$

and compactly as

$$\sum_{\ell=0}^{\nu-\kappa} L_E L_F^{\nu-\kappa-\ell}(\beta_\ell H(x)) + \beta_{\nu-\kappa+1} K(x)$$
$$= L_E L_F^{\nu-\kappa} q(x) + \sum_{\ell=1}^{\nu-\kappa} \alpha_\ell L_E L_F^{\nu-\kappa-\ell} q(x), \kappa = 1, \cdots, \nu$$
$$\beta_0 K(x) = 0 \qquad (3.5)$$

The above are exactly condition b) of the Proposition.

Application: *Non-isothermal chemical reactor monitoring*

Liquid-phase oxidation reactions are notorious for being highly exothermic and for involving serious safety threats. One well-studied example is the reaction of N-methyl pyridine (A) with hydrogen peroxide (B) in the presence of a catalyst. The product of the reaction, Methyl Pyridine N-Oxide is an important intermediate in several reactions in pharmaceutical industry including the production of anti-ulcer and anti-inflammatory drugs [16]. The key issue in the operation of liquid-phase oxidation reactors is safety: both the organic reactant is usually hazardous at high temperatures, and also hydrogen peroxide decomposition could pose serious safety threats.

It is assumed the reactor is well-mixed and has constant volume with an inlet stream containing N-methyl pyridine (A) + catalyst Z (assumed to be fully dissolved) and hydrogen peroxide (B). The catalyst is assumed to be completely dissolved in the pyridine stream and its concentration is assumed to be constant in the reactor The dynamics of the reactor [16] is described by:

$$\frac{dc_A}{dt} = \frac{F}{V}(c_{A_{in}} - c_A) - R(c_A, c_B, \theta, w_1(t)) \qquad (3.9)$$
$$\frac{dc_B}{dt} = \frac{F}{V}(c_{B_{in}} - c_B) - R(c_A, c_B, \theta, w_1(t))$$
$$\frac{d\theta}{dt} = \frac{F}{V}(\theta_{in} - \theta) + \frac{(-\Delta H)_R}{\rho c_p} R(c_A, c_B, \theta)$$
$$- \frac{(U + w_2(t))A}{\rho c_p V}(\theta - \theta_J)$$
$$\frac{d\theta_J}{dt} = \frac{F_J}{V_J}(\theta_{J_{in}} - \theta_J) + \frac{(U + w_2(t))A}{\rho_J c_{p_J} V_J}(\theta - \theta_J)$$

where $c_A$ and $c_B$ are the concentrations of Pyridine and Hydrogen Peroxide respectively in the reacting mixture, $\theta$ and $\theta_J$ are the temperatures of the reacting mixture and the jacket fluid respectively; these are the system states. The reaction rate $R(c_A, c_B, \theta, w_1(t))$ is an empirically derived expression and is given by $R(c_A, c_B, \theta, w_1(t)) = \frac{A_1 e^{-\frac{E_1}{\theta}} A_2 e^{-\frac{E_2}{\theta}} c_A c_B Z}{1 + A_2 e^{-\frac{E_2}{\theta}} c_B} + A_3 e^{-\frac{E_3}{\theta}} c_A c_B + w_1(t)$, in which $w_1(t)$ represents the uncertainty in the empirical expression. $c_{A,in}$ and $c_{B,in}$ are the feed concentrations of pyridine and hydrogen peroxide respectively, F and $F_J$ are the reactant feed and coolant flowrates respectively, $\theta_{in}$ and $\theta_{J,in}$ are the inlet temperatures for the reactor and the cooling jacket respectively, V and $V_J$ are the reactor volume and cooling jacket volume respectively $(-\Delta H)_R$ is the heat of reaction, $\rho, c_p$ and $\rho_J, c_{p_J}$ are the densities and heat capacities of reaction mixture and cooling fluid respectively, U and A are the overall heat transfer coefficient and heat transfer area respectively, and $w_2(t)$ represents the uncertainty in the heat transfer coefficient. $A_1, A_2, A_3$ and $E_1, E_2, E_3$ are the reaction rate parameters, pre-exponential factors and rescaled activation energies respectively, and Z is the catalyst concentration (constant).

When reactants are potentially hazardous, special precautions are taken in terms of using relatively dilute feeds and the reaction taking place at a relatively low temperature. In terms of monitoring the operation of the reactor, the temperature $\theta$ of the reacting mixture as well as the total sum of hazardous chemicals' concentrations $c_A + c_B$ are critical quantities to be monitored. Temperature is easy and inexpensive to measure, but concentrations generally need to be estimated from temperature measurements. Consider therefore the problem of building an observer for the dynamic system (3.9), driven by the temperature measurements

$$y_1 = \theta \qquad (3.10)$$
$$y_2 = \theta_J$$

the objective being to estimate the sum of the reactant concentrations

$$z = c_A + c_B \quad (3.11)$$

decoupled from disturbances $w_1$ and $w_2$. To derive the functional observer, it is convenient to perform appropriate translation of axes to shift the equilibrium point to the origin by defining $c'_A = c_A - c_{A,s}$, $c'_B = c_B - c_{B,s}$, $\theta' = \theta - \theta_s$, $\theta'_J = \theta_J - \theta_{J,s}$, where $(c_{A,s}, c_{B,s}, \theta_s, \theta_{J,s})$ is the steady state (equilibrium point) of the reactor.

$$\frac{dc'_A}{dt} = -\frac{F}{V}c'_A - [R(c'_A + c_{A,s}, c'_B + c_{B,s}, \theta' + \theta_s) - R(c_{A,s}, c_{B,s}, \theta_s, w_1)] \quad (3.12)$$

$$\frac{dc'_B}{dt} = -\frac{F}{V}c'_B - [R(c'_A + c_{A,s}, c'_B + c_{B,s}, \theta' + \theta_s) - R(c_{A,s}, c_{B,s}, \theta_s, w_1)]$$

$$\frac{d\theta'}{dt} = -\frac{F}{V}\theta' + \frac{(-\Delta H)_R}{\rho c_p}[R(c'_A + c_{A,s}, c'_B + c_{B,s}, \theta' + \theta_s) - R(c_{A,s}, c_{B,s}, \theta_s, w_1)] - \frac{(U + w_2(t))A}{\rho c_p V}(\theta' - \theta'_J)$$

$$\frac{d\theta'_J}{dt} = -\frac{F_J}{V_J}\theta'_J + \frac{(U+w_2(t))A}{\rho_J c_{p_J} V_J}(\theta' - \theta'_J)$$

For the above system, a scalar disturbance decoupled functional observer can be built ($\nu=1$), with the necessary and sufficient conditions (Proposition 2) being satisfied for the following choices of $\beta_0, \beta_1 \in \mathbb{R}^2, \alpha_1 \in \mathbb{R}$:

$$\beta_0 = \left[-\frac{2\rho c_p}{(-\Delta H)_R}, -\frac{2\rho_J c_{p_J} V_J}{(-\Delta H)_R V}\right],$$

$$\beta_1 = \left[-\frac{2\rho c_p}{(-\Delta H)_R}\frac{F}{V}, -\frac{2\rho_J c_{p_J} F_J}{(-\Delta H)_R V}\right], \alpha_1 = \frac{F}{V}$$

all the conditions are satisfied, leading to

$$T(c'_A, c'_B, \theta', \theta'_J) = c'_A + c'_B + \frac{2\rho c_p}{(-\Delta H)_R}\theta' + \frac{2\rho_J c_{p_J} V_J}{(-\Delta H)_R V}\theta'_J$$

and the functional observer:

$$\frac{d\hat{\xi}}{dt} = -\frac{F}{V}\hat{\xi} - \frac{2\rho_J c_{p_J} V_J}{(-\Delta H)_R V}\left(\frac{F_J}{V_J} - \frac{F}{V}\right)y'_2 \quad (3.13)$$

$$\hat{z} = \hat{\xi} - \frac{2\rho c_p}{(-\Delta H)_R}y'_1 - \frac{2\rho_J c_{p_J} V_J}{(-\Delta H)_R V}y'_2$$

For the following parameter values (see [16]):
$C_{A,in} = 4 \frac{mol}{l}$, $C_{B,in} = 3 \frac{mol}{l}$, $\theta_{in} = 333$ K, $\theta_{J,in} = 300$ K,
$F = 0.02 \frac{l}{min}$, $F_j = 1 \frac{l}{min}$, $V = 1$ l, $V_J = 3 \times 10^{-2}$ l,
$A_1 = e^{8.08}$ l mol$^{-1}$s$^{-1}$, $A_2 = e^{28.12}$ l mol$^{-1}$s$^{-1}$, $A_3 = e^{25.12}$ l mol$^{-1}$, $E_1 = 3952$ K, $E_2 = 7927$ K, $E_3 = 12989$ K,
$\Delta H_R = -160 \frac{kJ}{mol}$, $\rho = 1200 \frac{g}{l}$, $\rho_J = 1200 \frac{g}{l}$, $c_{p_J} = 3.4 \frac{J}{gK}$.
$c_p = 3.4 \frac{J}{gK}$, $UA = 0.942 \frac{W}{K}$, $Z = 0.0021 \frac{mol}{l}$
the corresponding reactor steady state is:

$c_{A,s} = 1.211 \frac{mol}{l}$, $c_{B,s} = 0.211 \frac{mol}{l}$, $\theta_s = 386.20$ K, $\theta_{J,s} = 300.02$ K,

and we have simulated the reactor start-up, under the following initial conditions:
$c_A(0) = 0, c_B(0) = 0, \theta(0) = 300$ K, $\theta_J(0) = 300$ K.
The system and the observer (3.12) and 3.13) were simulated using constant disturbances $w_1(t) = 10^{-4} \frac{mol}{l \cdot s}$, $w_2(t) = 0.2 \frac{W}{m^2 \cdot K}$. Figure 1 compares the functional observer's estimate $\hat{c}_A + \hat{c}_B = \hat{z} + c_{A,s} + c_{B,s}$ to the system's total reactant concentration $c_A + c_B$, and provides a plot of the corresponding estimation error, when the initialization error is $1 \frac{mol}{l}$.

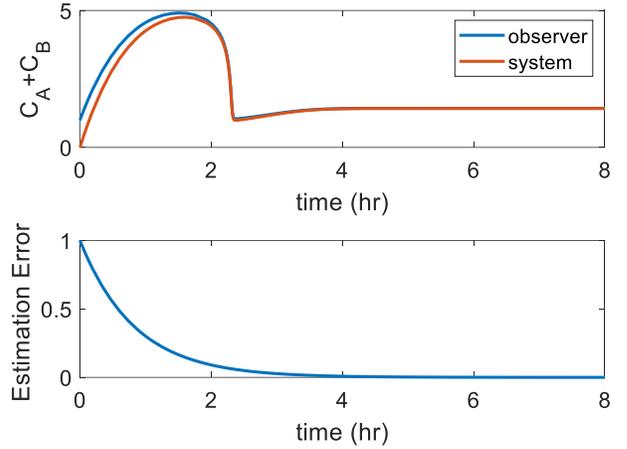

Figure 1: (a) System's and observer's response for $\hat{\xi}(0) - T(c'_A(0), c'_B(0), \theta'(0), \theta'_J(0)) = 1 \frac{mol}{\ell}$ (b) Estimation error $\hat{z}(t) - z(t)$

## IV. APPLICATION OF THE FUNCTIONAL OBSERVER TO THE PROBLEM OF SIMULTANEOUS FAULT DETECTION AND FAULT ESTIMATION

Functional observers of the form (2.2) are amenable to fault diagnosis applications. In this section, the problem of simultaneous disturbance decoupled detection and estimation of faults in nonlinear systems is considered. Fault is an unexcepted major deviation in process variables from the usual conditions. They can be categorized into different types based on their location. These include: (i) sensor faults (ii) component faults and (iii) actuator faults. Sensor faults may degrade performance of decision-making systems, including feedback control system, safety control system, quality control system, state estimation system, optimization system. The most common sensor faults include a) bias b) drift c) performance degradation d) sensor freezing e) calibration error A fault in an actuator may cause loss of control in automated control systems. Actuator faults include, for example, stuck-up of control valves and faults in pumps, etc. Several common faults in servomotors include Lock-in-Place

(LIP), Float, Hard-over Failure (HOF) and Loss of Effectiveness (LOE). Component faults occur in the equipment of plant. Examples for this could include leaks/blockages pipeline and tanks. These faults change the physical characteristics of the component and as a result can lead to significant change in the dynamics of the process.

It is assumed that the fault of interest originates from an exo-system of the following form

$$\frac{dx_o}{dt} = Rx_o \quad (4.1)$$
$$f = Qx_o$$

where $x_o \in \mathbb{R}^{n_o}, R \in \mathbb{R}^{n_o \times n_o}, Q \in \mathbb{R}^{1 \times n_o}$. Examples of faults that can be represented by the exo-system (4.1) include, step ($R = 0$, $Q = 1$), ramp ($R = \begin{bmatrix} 0 & 1 \\ 0 & 0 \end{bmatrix}, Q = [1\ 0]$) and sine wave $R = \begin{bmatrix} 0 & \omega \\ -\omega & 0 \end{bmatrix}, Q = [1\ 0]$. This approach is in the same vein as the one taken in the internal model principle in regulation problems to eliminate tracking errors in the presence of disturbances[13-15]. Equation (4.1) can be thought of as a disturbance generator[15].

Now consider a nonlinear system of the form

$$\frac{dx}{dt} = F(x) + G(x)f + E(x)W \quad (4.2)$$
$$y = H(x) + J(x)f + K(x)W$$
$$z = f$$

and where $f \in \mathbb{R}$ represents a potential fault arising from equipment malfunction, such that f is zero under normal operation, but f assumes a nonzero value in the event of a sudden malfunction. The quantity to be estimated (z) is now equal to the fault of interest.

The overall system of process (4.2) and exo-system (4.1) in cascade is

$$\frac{dx}{dt} = F(x) + G(x)Qx_o + E(x)W \quad (4.3)$$
$$\frac{dx_o}{dt} = Rx_o$$
$$y = H(x) + Qx_oJ(x) + K(x)W$$
$$z = Qx_o$$

To be able to detect the occurrence of a fault and at the same time estimate its size, a linear functional observer of the form

$$\frac{d\hat{\xi}}{dt} = A\hat{\xi} + By \quad (4.4)$$
$$\hat{f} = C\hat{\xi} + Dy$$

will be designed based on system (4.3) so that the response of the estimate $\hat{f}$ of the series connection of (4.3) and (4.4) will have the following properties
- it asymptotically approaches $f = Qx_o$
- it is unaffected by the disturbances W

Such an observer will exist if there exists a differentiable mapping $T: \mathbb{R}^{n+n_o} \to \mathbb{R}^\nu$ such that:

$$\frac{\partial T}{\partial x}(x, x_o)(F + G(x)Qx_o) + \frac{\partial T}{\partial x_o}(x, x_o)Rx_o$$
$$= AT(x, x_o) + B(H(x) + J(x)Qx_o) \quad (4.5)$$
$$CT(x, x_o) + D\big((H(x) + J(x)Qx_o)\big) = Qx_o \quad (4.6)$$

and in addition, it satisfies the disturbance decoupling conditions:

$$-\frac{\partial T}{\partial x}(x, x_o)E(x) + BK(x) = 0 \quad (4.7)$$

$$DK(x) = 0 \quad (4.8)$$

Propositions 3 and 4 that follow provide necessary and sufficient conditions for the functional observer (4.4) to satisfy the functional observer conditions (4.5) and (4.6), and the disturbance decoupling conditions (4.7) and (4.8), respectively. These are corollaries of Propositions 1 and 2 given in sections II and III.

Proposition 3: *There exists a functional observer of the form (4.4) for the system (4.3), with the eigenvalues of A being the roots of the polynomial $\lambda^\nu + \alpha_1 \lambda^{\nu-1} + \cdots + \alpha_{\nu-1}\lambda + \alpha_\nu$, if and only if there exist constant row vectors $\beta_0, \beta_1, \cdots, \beta_\nu \in \mathbb{R}^p$ that satisfy:*

$$L_{F_e}^\nu(\beta_0(H(x) + J(x)Qx_o)) + L_{F_e}^{\nu-1}(\beta_1(H(x) + J(x)Qx_o)) + \cdots$$
$$+ L_{F_e}(\beta_{\nu-1}(H(x) + J(x)Qx_o)) + \beta_\nu(H(x) + J(x)Qx_o) =$$
$$Q(R^\nu + \alpha_1 R^{\nu-1} + \cdots + \alpha_{\nu-1}R + \alpha_\nu I)x_o \quad (4.9)$$

where $L_{F_e} = \sum_{i=1}^n F_k(x)\frac{\partial}{\partial x_k} + Qx_o \sum_{i=1}^n G_k(x)\frac{\partial}{\partial x_k} + \sum_{i=1}^{n_o} R_k x_o \frac{\partial}{\partial x_{o_k}}$

It should be noted that under the conditions of Proposition 3, the mapping $T(x, x_o) = \begin{bmatrix} T_1(x, x_o) \\ \vdots \\ T_{\nu-1}(x, x_o) \\ T_\nu(x, x_o) \end{bmatrix}$, where:

$$T_1(x, x_o)$$
$$= Q(R^{\nu-1} + \alpha_1 R^{\nu-2} + \cdots + \alpha_{\nu-1}I)x_o$$
$$- L_{F_e}^{\nu-1}(\beta_0(H(x) + J(x)Qx_o)) - \cdots - \beta_{\nu-1}(H(x) + J(x)Qx_o)$$
$$\vdots \quad (4.10)$$
$$T_{\nu-1}(x, x_o)$$
$$= Q(R + \alpha_1 I)x_o - L_{F_e}(\beta_0(H(x) + J(x)Qx_o))$$
$$- \beta_1(H(x) + J(x)Qx_o)$$
$$T_\nu(x, x_o) = Qx_o - \beta_0(H(x) + J(x)Qx_o)$$

satisfies conditions (4.5) and (4.6) with A, B, C, D given by (2.8).

The results of Proposition 3 can now be specialized for some common fault types.

→ For the special case of a step fault, the exo-system (4.1) is $\frac{dx_o}{dt} = 0$ and $f = x_o$ where $x_o$ is scalar, and condition (4.9) takes the form:

$$L_{F_e}^\nu(\beta_0(H(x) + J(x)x_o)) + \cdots + L_{F_e}(\beta_{\nu-1}(H(x) + J(x)x_o))$$
$$+ \beta_\nu(H(x) + J(x)x_o) = \alpha_\nu x_o$$

It should be noted that in this case, the above condition can be rescaled so as to be independent of the eigenvalue-dependent coefficient $\alpha_\nu$:

$$L_{F_e}^\nu(v_0(H(x) + J(x)x_o)) + \cdots + L_{F_e}(v_{\nu-1}(H(x) + J(x)x_o))$$
$$+ v_\nu(H(x) + J(x)x_o) = x_o$$

where $v_j = \frac{\beta_j}{\alpha_\nu}$ for $j = 1,2..,\nu$.

This implies that the eigenvalues can be freely assigned to ensure stability of the error dynamics.

→ For the special case of a ramp, R and Q in system (4.1) are $\begin{bmatrix} 0 & 1 \\ 0 & 0 \end{bmatrix}$ and $[1\ 0]$ respectively, and (4.9) becomes:

$$L_{F_e}^\nu(\beta_0(H(x) + J(x)x_{o1})) + \cdots + L_{F_e}(\beta_{\nu-1}(H(x) + J(x)x_{o1})) + \beta_\nu(H(x) + J(x)x_{o1}) = \alpha_{\nu-1}x_{o2} + \alpha_\nu x_{o1}$$

with the understanding that for $\nu = 1$, $\alpha_0 = 1$.

Here, rescaling of the row vectors $\beta_0, \ldots, \beta_\nu$ to make the condition independent of eigenvalues is not possible.

Step type faults are a special case of ramp faults in which $x_{o2} = 0 \; \forall \; t$. Therefore, a functional observer built to detect ramp faults can also detect step faults.

<u>Proposition 4</u>: *Suppose that there exist constant row vectors $\beta_0, \beta_1, \cdots, \beta_\nu \in \mathbb{R}^p$ that satisfy condition (4.9) and that the matrices (A, B, C, D) have been chosen according to (2.8), so that (4.5) and (4.6) hold with T(x) given by (4.10). The functional observer (4.4) will satisfy the disturbance decoupling conditions (4.7) and (4.8) if and only if*

$$\sum_{\ell=0}^{\nu-\kappa} L_E L_{F_e}^{\nu-\kappa-\ell}(\beta_\ell [(H(x) + J(x)x_o)]) + \beta_{\nu-\kappa+1} K(x) = 0,$$
$$\kappa = 1, \cdots, \nu \qquad (4.11)$$
$$\beta_0 K(x) = 0$$

The disturbance decoupled fault detection and estimation approach developed in this section can be directly extended to the case of multiple faults, so that in addition, fault isolation is accomplished. For example, in a system with $n_f$ faults, $n_f$ functional observers can be designed, one for each fault, where for functional observer i, fault i is to be detected and all the other $n_f - 1$ faults are disturbances that are decoupled (see Figure 2).

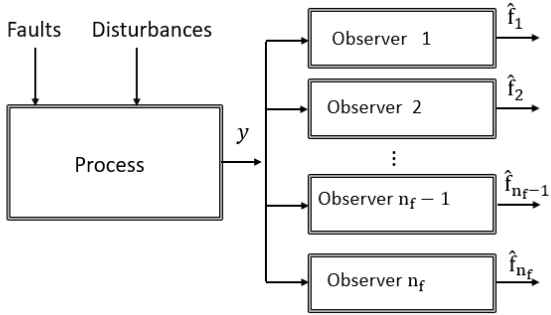

Figure 2 Fault isolation and estimation scheme, based on a set of observers, one for each fault

<u>Application</u>: *Fault detection and estimation for a chemical reactor*

Consider the CSTR reactor presented before with an added concentration measurement and possible persistent faults in the inlet temperature of the coolant and the concentration sensor. The uncertainties in the pre-exponential factors of the reaction rate are same as before while no uncertainty in the heat transfer coefficient is assumed. The dynamics are as follows

$$\frac{dc'_A}{dt} = -\frac{F}{V} c'_A - [R(c'_A + c_{A,s}, c'_B + c_{B,s}, \theta' + \theta_s) - R(c_{A,s}, c_{B,s}, \theta_s, w_1)] \quad (4.12)$$

$$\frac{dc'_B}{dt} = -\frac{F}{V} c'_B - [R(c'_A + c_{A,s}, c'_B + c_{B,s}, \theta' + \theta_s) - R(c_{A,s}, c_{B,s}, \theta_s, w_1)]$$

$$\frac{d\theta'}{dt} = -\frac{F}{V} \theta' + \frac{(-\Delta H)_R}{\rho c_p} [R(c'_A + c_{A,s}, c'_B + c_{B,s}, \theta' + \theta_s) - R(c_{A,s}, c_{B,s}, \theta_s, w_1)] - \frac{UA}{\rho c_p V} (\theta' - \theta'_J)$$

$$\frac{d\theta'_J}{dt} = -\frac{F_J}{V_J} \left(\theta'_J + f_2(t)\right) + \frac{UA}{\rho_J c_{pJ} V_J} (\theta' - \theta'_J)$$

$$y'_1 = C'_A + f_1(t)$$
$$y'_2 = \theta'$$
$$y'_3 = \theta'_J$$

where $f_2$ is a possible step-fault that originates from a potential malfunction of the coolant feeding system and $f_1$ represents a potential ramp type fault in the analytical sensor. Like before, the model is converted to deviation form. The goal is to design a fault diagnosis scheme that can detect, isolate and estimate these two persistent step faults $f_1$ and $f_2$ in the presence of uncertainties in the reaction rate. To this end, two scalar functional observers are built (i) to estimate the analytical sensor fault ($f_1$) while considering $f_2$ as an additional disturbance. (ii) to estimate inlet jacket temperature fault $f_2$ considering $f_1$ as an additional disturbance.

<u>Functional Observer 1</u>: Estimation of the analytical sensor fault $f_1$ while considering $f_2$ as an additional disturbance. In this case, the functional observer will be based on an extended system by appending the equation $\frac{dx_{o1}}{dt} = x_{o2}, \frac{dx_{o2}}{dt} = 0, f_1 = x_{o1}$ to system (4.12).

A scalar ($\nu=1$) functional observer can be designed with the following choice of $\beta_0$ and $\beta_1$:

$$\beta_1 = \left[\alpha_1, \frac{\left(\frac{F}{V} + \frac{UA}{\rho c_p V}\right) \rho c_p}{-\Delta H_R}, \frac{-UA}{(-\Delta H_R)V}\right] \quad (4.13)$$

$$\beta_0 = \left[1, \frac{\rho c_p}{-\Delta H_R}, 0\right]$$

and design parameters $A = -\alpha_1 = -\frac{F}{V}$, $B = \beta_1 - \alpha_1 \beta_0$, $C = 1$, $D = \beta_0$:

$$\frac{d\hat{\xi}}{dt} = -\frac{F}{V} \hat{\xi} + \left(\frac{UA}{(-\Delta H_R)V}\right) y'_2 - \frac{UA}{(-\Delta H_R)V} y'_3$$

$$\hat{f}_1 = \hat{\xi} + y'_1 + \frac{\rho c_p}{(-\Delta H_R)} y'_2 \quad (4.14)$$

<u>Functional Observer 2</u>: Estimation of the inlet cooling jacket temperature fault $f_2$ considering $f_1$ as an additional disturbance. In this case, the functional observer will be based on an extended system by appending the equation $\frac{df_2}{dt} = 0$ to system (4.12)

A scalar ($\nu=1$) functional observer can be designed with the following choice of $\beta_0$ and $\beta_1$:

$$\beta_1 = \alpha_1 \left[0, -\frac{UA}{\rho_J C_{pJ} F_J}, 1 + \frac{UA}{\rho_J C_{pJ} F_J}\right] \quad (4.15)$$

$$\beta_0 = \alpha_1 \left[0, 0, \frac{V_J}{F_J}\right]$$

and design parameters $A = -\alpha_1 = -0.01$, $B = \beta_1 - \alpha_1 \beta_0$, $C = 1$, $D = \beta_0$

$$\frac{d\hat{\xi}}{dt} = -0.01\hat{\xi} - \left(0.01\frac{UA}{\rho_j c_{p_j} F_j}\right)y_2'$$
$$- \left(0.01^2 \frac{V_j}{F_j} - 0.01 - 0.01\frac{UA}{\rho_j c_{pj} F_j}\right)y_3'$$
$$\hat{f}_2 = \hat{\xi} + 0.01\frac{V_j}{F_j}y_3' \qquad (4.16)$$

The two functional observers are simulated for the following scenario:
Two persistent faults
$$f_1(t) = \begin{cases} 0, & t < 2000s \\ 0.001t - 1, & t \geq 2000s \end{cases}, f_2(t) = \begin{cases} 0, & t < 5000s \\ 10, & t \geq 5000s \end{cases}$$
are assumed to occur along with the following disturbances $w_1(t) = 10^5, w_2(t) = 10^5$. It is assumed that initially the process is operating at steady state. An initialization error $(\hat{\xi}(0)-T(x(0))= 1$ is assumed for both observers (where T(x) is given by equation 4.10). The fault estimates are plotted in Figure 3.

Both estimates from time t =0 to 2000s are decay to 0 in the absence of faults. When the sensor fault occurs at time t=2000s a deviation is seen in $\hat{f}_1$ whereas $\hat{f}_2$ is equal to 0. At time t=5000s a deviation is observed in $\hat{f}_2$ indicating the presence of a fault in thse inlet coolant temperature. Both profiles eventually converge to their actual fault profiles.

It is to be noted that the step faults represent special cases of ramp faults (with slope =0) and hence, observers designed to detect ramp faults can also detect step faults.

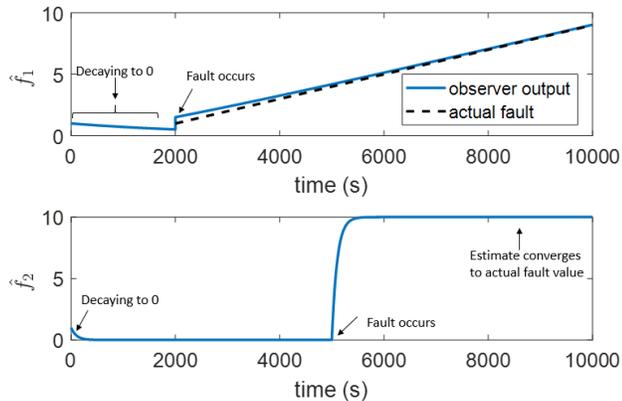

Figure 3: Fault estimates vs time with initialization error=1. The fault estimates converge to their actual trajectories

## V. CONCLUSIONS

This work has studied the design of disturbance decoupled functional observers from the point of view of observer error linearization. Necessary and sufficient conditions are derived for disturbance decoupling of the functional observer output. A fault estimation methodology based on disturbance decoupled functional observers is proposed. Throughout the study, the methods are applied to a non-isothermal CSTR case study.